\documentclass[10pt,leqno]{article}

\usepackage{amsmath,amsfonts,amscd,amssymb,theorem}

\long\def\comment#1\endcomment{}

\makeatletter
\begingroup
\gdef\th@dotted{\normalfont\itshape
  \def\@begintheorem##1##2{%
        \item[\hskip\labelsep \theorem@headerfont ##1\ ##2.]}%
\def\@opargbegintheorem##1##2##3{%
   \item[\hskip\labelsep \theorem@headerfont ##1\ ##2\ (##3).]}}
\endgroup
\makeatother

\theoremstyle{dotted}

\newtheorem{theorem}{Theorem}[section]
\newtheorem{lemma}[theorem]{Lemma}
\newtheorem{conj}[theorem]{Conjecture}

\makeatletter
\begingroup
\gdef\th@upshape{\normalfont
  \def\@begintheorem##1##2{%
        \item[\hskip\labelsep \theorem@headerfont ##1\ ##2.]}%
\def\@opargbegintheorem##1##2##3{%
   \item[\hskip\labelsep \theorem@headerfont ##1\ ##2\ (##3).]}}
\endgroup
\makeatother

\theoremstyle{upshape}

\newtheorem{defn}[theorem]{Definition}
\newtheorem{remark}[theorem]{Remark}

\makeatletter
\renewcommand{\subsection}{\@startsection{subsection}{2}{0pt}{-3ex
plus -1ex minus -0.2ex}{-2mm plus -0pt minus
-2pt}{\normalfont\bfseries}} \makeatother

\makeatletter
\@addtoreset{equation}{section}
\makeatother

\newcommand{\cntrct}                
{\hspace{2pt}\raisebox{1pt}{\text{$\lrcorner$}}\hspace{2pt}}

\newcommand{\proof}[1][Proof.]{\smallskip\noindent{\em #1}}
\def\endproof{\hfill\ensuremath{\square}\par\medskip}

\def\eqref#1{\thetag{\ref{#1}}}

\let\latexref=\ref
\def\ref#1{{\normalfont{\latexref{#1}}}}

\newcommand{\idot}{{\:\raisebox{1pt}{\text{\circle*{1.5}}}}}
\newcommand{\hdot}{{\:\raisebox{3pt}{\text{\circle*{1.5}}}}}
\newcommand{\calo}{{\cal O}}

\newcommand{\gr}{\operatorname{\sf gr}}

\newcommand{\C}{{\mathbb C}}
\newcommand{\R}{{\mathbb R}}
\newcommand{\A}{{\mathbb A}}
\newcommand{\Z}{{\mathbb Z}}

\newcommand{\E}{{\mathcal E}}
\newcommand{\T}{{\mathcal T}}
\newcommand{\N}{{\mathcal N}}

\newcommand{\X}{{\mathfrak X}}
\newcommand{\g}{{\mathfrak g}}

\newcommand{\rk}{\operatorname{\sf rk}}
\renewcommand{\dim}{\operatorname{\sf dim}}
\newcommand{\codim}{\operatorname{\sf codim}}
\newcommand{\cchar}{\operatorname{\sf char}}

\title{Sommese Vanishing for non-compact manifolds}

\author{D. Kaledin\thanks{Partially supported by CRDF grant
RM1-2354-MO02.}}

\begin{document}

\maketitle


\section*{Introduction}

Let $X$ be a smooth algebraic manifold of dimension $\dim X = n$
over an algebraically closed field $k$ of characteristic $\cchar k =
0$. Assume that $X$ is equipped with a birational projective map
$\pi:X \to Y$ to a normal irreducible algebraic manifold $Y$. The
Kodaira Vanishing Theorem admits a well-known generalization to this
relative situation, namely, the Grauert-Riemenschneider Vanishing
Theorem, which claims that the higher direct image sheaves
$R^k\pi_*\Omega_X^n$ on $Y$ are trivial for $k > 0$. One could
conjecture that the stronger Kodaira-Nakano Theorem can also be
generalized to the relative situation, so that one would have
$R^p\pi_*\Omega_X^q=0$ for $p+q > n$. However, a moment's reflection
-- take $X$ to be the blowup of a smooth point -- shows that this is
not true, unless one imposes some additonal assumptions on $X$ or on
the map $\pi$.

The goal of this paper is to prove the relative version of
Kodaira-Nakano Theorem for one set of such additional
assumptions. We will prove (see Theorem~\ref{main} below) that
$R^p\pi_*\Omega_X^q=0$ if $p+q > \dim X \times_Y X$, the dimension
of the fibered product of $X$ with itself over $Y$. We do not need
the map $X \to Y$ to be birational, projective is enough. In the
case when $\dim X \times_Y X = X$, the condition on degrees is the
original one, that is, $p + q > n$. Such maps $X \to Y$ are
necessarily generically finite, and they are known as {\em semismall
maps}.

A vanishing theorem of the same type has been proved some time ago
by A. Sommese \cite{S}, but only in the case when the base manifold
$Y$ is compact. Sommese Theorem is global -- rather than consider
the direct images, he considers the global cohomology groups
$H^p(X,\Omega_X^q \otimes \pi^*M)$, where $M$ is an ample line
bundle on $Y$. There is a well-known procedure for deducing the
vanishing of higher direct images from such a global vanishing, but
there is a hitch -- in order to apply it, one has to compactify $X$,
$Y$ and $\pi:X \to Y$ so that all assumptions on $\pi$ are
preserved. Unfortunately, usually it is not possible to compactify a
semismall map so that it stays semismall.

Our proof uses essentially the same ideas, but we arrange the
induction on dimension in a slightly different way; this allows us
to use an arbitrary smooth compactification but make all the
unpleasantness stay at infinity, where it belongs. The
characteristic $0$ assumption is used in two places: firstly, we
need it to use Hodge Theory, secondly, we need the resolution of
singularities to construct smooth compactifications. The Hodge
Theory part is not critical, since it can be done in characteristic
$p$ by the famous method of Delign-Illusie \cite{DI} (at the cost of
some additional assumptions on $X$ such as liftability to the ring
of second Witt vectors). But the resolution of singularities seems
essential.

As an application, we prove a general theorem on the topology of
symplectic contractions over $\C$ (Theorem~\ref{sympl}). Roughly
speaking, we assume that $Y$ is affine, $\pi:X \to Y$ is birational,
and $X$ carries a non-degenerate closed $2$-form $\Omega$, and we
deduce that for the fiber $E_y = \pi^{-1}(y) \subset X$ of the map
$\pi:X \to Y$ over any closed point $y \in Y$, the cohomology group
$H^k(E_y,\C)$ is trivial if $k=2p+1$, and carries a pure Hodge
structure of type $(p,p)$ if $k=2p$. We note that in the particular
case of the so-called {\em Springer resolution}, this statement has
been proved by a direct geometric argument in \cite{lu}.

\medskip

The author is definitely not an expert in the field of vanishing
theorems. For those like him, he would like to mention at this point
that there is an excellent and standard reference for all things
related to algebraic vanishing theorems, namely, the book \cite{EV}
by H. Esnault and E. Viehweg. In particular, it is from this book
that the author has learned the Sommese Theorem, as well as its
proof and all the ideas needed for the proof. The present paper has
also been very much influenced by the beautiful paper \cite{dCM}
devoted to the geometry of semismall maps.

The fibered product $X \times_Y X$ appears in the statements for
purely numerical reasons. It would be interesting to see if it has
any real singificance.

\subsection*{Acknowledgements.} This paper is a part of an ongoing
joint project with R. Bezrukavnikov on algebraic symplectic
manifolds and their quantization; although he didn't take part in
this particular effort, it is easily understood that his implicit
contribution goes way beyond the usual friendly help. I am very much
indebted to M. de Cataldo for explaining to me the results of
\cite{dCM} and for continuously answering my questions, some of them
offensively stupid. Part of the work was done while I was visiting
the University of Mainz at the invitation of M. Lehn. It is a
pleasure to thank this wonderful institution. It is even greater
pleasure to thank M. Lehn for his hospitality, patience and
expertise. Many thanks are due to E. Amerik, R. Bezrukavnikov and
A. Kuznetsov, who have read a preliminary version of this paper and
suggested several improvements. Finally, I am grateful to M. de
Cataldo for indicating a mistake in the proof of Lemma~\ref{redu},
and to V. Ginzburg for attracting my attention to the paper
\cite{lu}.

\section{Statements and preliminaries.}

Fix an algebraically closed field $k$ of characteristic $\cchar k =
0$. Let $X$ be a regular scheme over $k$ of dimension $\dim X =
n$. Assume given a projective map $\pi:X \to Y$ from $X$ to a normal
irreducible algebraic variety $Y$ over $k$. Denote by $X \times_Y X$
the fibered product of $X$ with itself over $Y$. Here is our main
theorem.

\begin{theorem}\label{main}
For any $p,q \geq 0$ such that $p + q > \dim X \times_Y X$ we have
$$
R^p\pi_*\Omega^q_X = 0.
$$
\end{theorem}

Before we start proving it, we need to set up some preliminary
machinery. First of all, we explain the numerical bound in
Theorem~\ref{main}.

\begin{lemma}\label{cohodim}
Assume given a smooth variety $X$ over $k$ equipped with a proper
map $\pi:X \to Y$ into an affine variety $Y$. Then
$$
H^l_{DR}(X) = 0
$$
whenever $l > \dim X \times_Y X$.
\end{lemma}

\proof{} Since $k$ has characteristic $0$, we may assume $k = \C$
and work in the analytic topology. We can compute the cohomology
$H^\hdot_{DR}(X) = H^\hdot(X,\C)$ by applying the Leray spectral
sequence for the map $\pi:X \to Y$. Let $Y_p \subset Y$ be the
closed subvariety of points $y \in Y$ such that $\dim \pi^{-1}(y)
\geq p$. By proper base change, the sheaf $R^k\pi_*\C$ is supported
on $Y_p$, where $p$ is the smallest integer such that $k \leq 2p$.
Since $Y_p \subset Y$ are affine varieties, the group
$$
H^{l-k}(Y,R^k\pi_*\C) = H^{l-k}(Y_p, R^k\pi_*\C)
$$
vanishes whenever $l - k > \dim Y_p$. On the other hand, $\dim X
\times_Y X \geq \dim Y_p + 2p$, so that $l > \dim X \times_Y X$
implies $l > \dim Y_p + 2p$. Collecting all this together, we see
that
$$
H^{l-k}(Y,R^k\pi_*\C)=0
$$
whenever $l > \dim X \times_Y X$.
\endproof

\noindent
Next, assume given a smooth manifold $X$ over $k$ equipped with a
line bundle $L$.

\begin{defn}
A {\em flag} for $L$ is a sequence
$$
\emptyset = W_l \subset \dots \subset W_1 \subset W_0 = X
$$
such that $W_{i+1}$ is either empty or a Cartier divisor in $W_i$,
and the line bundle $\calo_{W_i}(W_{i+1})$ on $W_i$ is isomorphic to
the restriction of the line bundle $L$ onto $W_i$. A flag is called
{\em smooth} if all the $W_i$ are smooth.
\end{defn}

If in addition $X$ is equipped with a simple normal crossing divisor
$D$, then a flag $W_i$ is said to be {\em transversal to $D$} if
every $W_i$ intersects transversally with every intersection of the
components of the divisor $D$. If this is the case, then in
particular the union $D \cup W_1$ is a simple normal crossing
divisor in $X$. Moreover, by induction $(D \cap W_i) \cup W_{i+1}$
is a normal crossing divisor in $W_i$ for every $i$.

\begin{lemma}\label{flag}
Assume given a line bundle $L$ on a smooth manifold $X$ over $k$. If
$L$ is base-point free, it admits a smooth flag. Moreover, if $X$ is
equipped with a normal crossing divisor $D$, then such a flag can be
chosen so that it is transversal to $D$.
\end{lemma}

\proof{} By induction on $\dim X$, it suffices to find a smooth
divisor $W_1 \subset X$ in the linear system $|L|$, transversal to
$D$ if necessary. By the Bertini Theorem, a generic member of $|L|$
will do.
\endproof

\begin{remark} If the line bundle $L$ is very ample, a flag for $L$
is of length $n = \dim X$. However, this need not be so if $L$ is
only base-point free. In the extreme case of trivial $L$, a flag
consists of two subvarieties: $W_0 = X$ and $W_1 =
\emptyset$. Nevertheless, it is a smooth flag.
\end{remark}

\begin{lemma}\label{twtolog}
Let $X$ be a smooth manifold over $k$ equipped with a normal
crossing divisor $D$, a line bundle $L$ and a smooth flag $W_i$ for
$L$ which is transversal to $D$. Then the logarithmic de Rham
sheaves $\Omega^\hdot_X\langle D \rangle \otimes L$ admit an
increasing filtration $F_i$ such that
$$
\gr_F^i\left(\Omega^\hdot_X\langle D \rangle \otimes L\right) =
\Omega^\hdot_{W_i}\langle (D \cap W_i) \cup W_{i+1} \rangle.
$$
\end{lemma}

\proof{} By induction on $\dim X$, it suffices to prove that the
natural restriction map $\Omega^\hdot_X \otimes L \to
\Omega^\hdot_{W_1} \otimes L$ and the embedding $\Omega^\hdot\langle
D \cup W_1 \rangle \subset \Omega^\hdot\langle D \rangle \otimes L$
induced by the embedding $\calo \subset \calo(W_1) \cong L$ together
give a short exact sequence
\begin{equation}\label{filt}
\def\longto{\ \to \ }
0 \longto \Omega^\hdot_X\langle D \cup W_0 \rangle \longto
\Omega^\hdot_X\langle D \rangle \otimes L \longto
\Omega^\hdot_{W_1}\langle D \cap W_1 \rangle \otimes L
\longto 0.
\end{equation}
This sequence is well-known, see e.g. \cite[Property 2.3c]{EV}. We
give a proof for the sake of completeness. The claim is local, so
that we may assume that the manifold $X$ is $\A^n$ with coordinates
$z_1,\dots,z_n$, that $W_1 \subset X$ is the coordinate hyperplane
$z_1 = 0$, and that $D \subset X$ is the union of the coordinate
hyperplanes $z_i = 0$, $i=2,\dots,l$ for some $l \leq n$. We have
$$
\Omega^\hdot_X \langle D \cup W_1 \rangle =
\Omega^\hdot_{\A^1}\langle o \rangle \boxtimes
\left(\Omega^\hdot_{\A^1}\langle o \rangle \boxtimes \dots \boxtimes
\Omega^\hdot_{\A^1}\langle o \rangle\right) \boxtimes
\left(\Omega^\hdot_{\A^1} \boxtimes \dots \boxtimes
\Omega^\hdot_{\A^1}\right),
$$
where $\Omega^\hdot_{\A^1}\langle o \rangle$ is the logarithmic de
Rham complex of $\A^1$ with singularities at the origin point $o \in
\A^1$, and the logarithmic factors correspond to coordinates $z_1$
and $z_2,\dots,z_l$. Moreover, we have
$$
\Omega^\hdot_X \langle D \rangle \otimes L =
\left(\Omega^\hdot_{\A^1} \otimes \calo(o)\right) \boxtimes
\left(\Omega^\hdot_{\A^1}\langle o \rangle \boxtimes \dots \boxtimes
\Omega^\hdot_{\A^1}\langle o \rangle\right) \boxtimes
\left(\Omega^\hdot_{\A^1} \boxtimes \dots \boxtimes
\Omega^\hdot_{\A^1}\right),
$$
and
$$
\Omega^\hdot_{W_1} \langle (D \cap W_1) \rangle \otimes L = k_o
\boxtimes \left(\Omega^\hdot_{\A^1}\langle o \rangle \boxtimes \dots
\boxtimes \Omega^\hdot_{\A^1}\langle o \rangle\right) \boxtimes
\left(\Omega^\hdot_{\A^1} \boxtimes \dots \boxtimes
\Omega^\hdot_{\A^1}\right),
$$
where $k_o$ is the skyscraper sheaf supported at the
origin. Therefore \eqref{filt} follows from the obvious exact
sequence
$$
\begin{CD}
0 @>>> \Omega^\hdot_{\A^1}\langle o \rangle @>>> \Omega^\hdot_{\A^1}
\otimes \calo(o) @>>> k_o @>>> 0.
\end{CD}
$$
This finishes the proof.
\endproof

\section{Proofs.}

We can now start proving Theorem~\ref{main}. First we reduce the
problem to a particular case.

\begin{lemma}\label{redu}
To prove Theorem~\ref{main}, it suffices to prove that for every
$p$, $q$ with $p + q > \dim X \times_Y X$, the sheaf
$$
R^p\pi_*\Omega^q_X
$$
vanishes if it is supported set-theoretically on a finite union $Y_0
\subset Y$ of closed points in $Y$.
\end{lemma}

\proof{} Apply induction on $n = \dim X$. Assume that
Theorem~\ref{main} is proved for all $X'$ with $\dim X' < n$, and
assume that support of the sheaf $R^p\pi_*\Omega^q_X$ has an
irreducible component $Y_0 \subset Y$ of dimension $r = \dim Y_0 >
0$. Shrinking $Y$ if necessary, we may assume that $Y_0$ is smooth
and affine, and that $Y$ admits a smooth projection $\rho_0:Y \to
\A^r$ into the affine space $\A^r$ which restricts to an \'etale map
$\rho _0:Y_0 \to U$ from $Y_0 \subset Y$ onto an open subset $U
\subset \A^r$. Shrinking $Y$ even further, we can also assume that
the associated projection $\rho = \rho_0 \circ \pi:X \to Y_0$ is a
smooth map. Moreover, we may assume that $R^p\pi_*\Omega^q_X$ is
supported on $Y_0 \subset Y$. The smooth map $\rho:X \to U$ induces
a filtration on $\Omega^q_X$ with associated graded pieces
\begin{equation}\label{rela}
\Omega^k(X/U) \otimes \rho^*\Omega^l_{U}, \qquad k+l =q,
\end{equation}
where $\Omega^k(X/U)$ is the sheaf of relative $k$-forms on
$X/U$. The projection formula gives
$$
R^p\pi_*\left(\Omega^k(X/U) \otimes \rho^*\Omega^l_U\right) =
R^p\pi_*\Omega^k_X \otimes \Omega^l_U.
$$
This sheaf tautologically vanishes whenever $l > r$. For any closed
point $u \in U$, the fiber $X_u = \rho^{-1}(u) \subset X$ is smooth
and equipped with a projective map $\pi:X_u \to Y_u$ into the fiber
$Y_u = \rho_0^{-1}(u) \subset Y$. By induction on $n = \dim X$, we
may assume that
$$
R^p\pi_*\Omega^k_{X_u} = 0
$$
whenever $p+k > \dim X_u \times_{Y_u} X_u$. By base change, this
implies that
$$
R^p\pi_*\Omega^k(X/U) = 0
$$ 
if $p+k > \dim X_u \times_{Y_u} X_u$ for all $u \in U$. By our
assumptions, for a generic $u \in U$ we have
$$
\dim X \times_Y X = \dim Y_0 + \dim X_u \times_{Y_u} X_u.
$$
Shrinking $U$ if necessary, we may assume that this holds for any
closed point $u \in U$. Therefore $p+q > \dim X \times_Y X$ and $l
\leq \dim Y_0$ implies $k+p = p+q -l > \dim X_u \times_{Y_u}
X_u$. We conclude that
$$
R^p\pi_*\left(\Omega^k(X/Y_0) \otimes \rho^*\Omega^l_U\right) = 0
$$
whenever $p+q = p+k+l > \dim X \times_Y X$. Applying the spectral
sequence associated to the filtration \eqref{rela}, we conclude that
indeed, $R^p\pi_*\Omega^q_X = 0$.
\endproof

\proof[Proof of Theorem~\ref{main}.] By Lemma~\ref{redu}, we may
assume that all the sheaves $R^p\pi_*\Omega^q_X$ with $p + q > \dim
X \times_Y X$ are supported in a finite subset in $Y$. Without loss
of generality we may also assume that $Y$ is affine. Choose a
projective variety $\overline{Y}$, a smooth projective variety
$\overline{X}$ and a projective map $\pi:\overline{X} \to
\overline{Y}$ so that $Y \subset \overline{Y}$ is a dense open
subset in $\overline{Y}$, $X = \pi^{-1}(Y) \subset \overline{X}$ is
a dense open subset in $\overline{X}$, $\pi:\overline{X} \to
\overline{Y}$ extends the given map $\pi:X \to Y$, and the
complement $D = \overline{X} \setminus X$ is a simple normal
crossing divisor in $\overline{X}$. Choose an ample line bundle $M$
on $\overline{Y}$. Let $l >> 0$ be an integer large enough so that
the sheaves
$$
R^p\pi_*\Omega^q_{\overline{X}} \otimes M^{\otimes l}
$$
are acyclic for all $p,q$. Replace $M$ with $M^{\otimes l}$.

The line bundle $L = \pi^*M$ on $\overline{X}$ is base-point free.
By Lemma~\ref{flag}, there exists a smooth flag $W_i$ for the bundle
$L$ which is transversal to $D$. By construction, we have $W_i =
\pi^{-1}(Z_i)$, where $Z_i$ form a flag for $M$ on
$\overline{Y}$. In particular, for any $i$ the complement $W_i
\setminus W_{i+1}$ comes equipped with a projective map
$$
\pi:W_i \setminus W_{i+1} \to Z_i \setminus Z_{i+1},
$$
and the intersection $\left(Z_i \setminus Z_{i+1}\right) \cap Y$ is
affine. By Lemma~\ref{cohodim}, we have
$$
H_{DR}^l\left(\left(W_i \setminus W_{i+1}\right) \cap X \right) = 0
$$
whenever $l > \dim W_i \times_{Z_i} W_i$. If $l > \dim X \times_Y
X$, then the vanishing holds for all $i$.

The cohomology $H_{DR}^l((W_i \setminus W_{i+1})\cap X)$ can be
computed by means of the logarithmic de Rham complex
$\Omega^\hdot_{W_i}\langle(D \cap W_i) \cup
W_{i+1}\rangle$. Moreover, the Hodge-de Rham spectral sequence for
this complex degenerates, so that we have
$$
H^p\left(\overline{X},\Omega^q_{W_i}\langle(D \cap W_i) \cup
W_{i+1}\rangle\right) = 0
$$
whenever $l = p+q > X \times_Y X$. Applying Lemma~\ref{twtolog}, we
deduce that
\begin{equation}\label{logvan}
H^p\left(\overline{X},\Omega_{\overline{X}}^q\langle D \rangle
\otimes L\right) = 0
\end{equation}
under the same assumption on $p$, $q$.

Consider now the sheaves $\Omega^\hdot_{\overline{X}} \otimes L$ on
$\overline{X}$ and the Leray spectral sequence for the map
$\pi:\overline{X} \to \overline{Y}$. By our assumption on $L =
\pi^*M$, the spectral sequence degenerates, and we have
$$
H^p\left(\overline{X},\Omega^q_{\overline{X}} \otimes L\right) =
H^0\left(\overline{Y},R^p\pi_*\Omega^q_{\overline{X}} \otimes M\right)
$$
for every $p$, $q$. Since $Y \subset \overline{Y}$ is affine, the
same degeneration holds over $Y$. 

Assume from now on that $p+q > \dim X \times_Y X$. By our general
assumption, the sheaf $R^p\pi_*\Omega^q_X$ is supported in a finite
union of closed points in $Y$. In particular, its support is
distinct from the complement $\overline{Y} \setminus Y$, and it is
therefore a direct summand in the sheaf
$R^p\pi_*\Omega^q_{\overline{X}}$ on $\overline{Y}$. We conclude
that the natural restriction map
\begin{multline}\label{rest}
H^p\left(\overline{X},\Omega^q_{\overline{X}} \otimes L\right) =
H^0(\overline{Y},R^p\pi_*\Omega^q_{\overline{X}}) \longrightarrow\\
\longrightarrow H^p\left(X,\Omega^q_X \otimes L\right) =
H^0(Y,R^p\pi_*\Omega^q_X)
\end{multline}
is surjective.

On the other hand, denote by $j:X \to \overline{X}$ the open
embedding. Since $D = \overline{X} \setminus X$ is a divisor in
$\overline{X}$, the map $j$ is affine, so that the higher direct
images $R^kj_*\E$ are trivial for $k \geq 1$ and arbitrary coherent
sheaf $\E$ on $X$. Therefore the restriction map \eqref{rest} is
induced by the embedding of sheaves $\Omega^q_{\overline{X}} \otimes
L \subset j_*\Omega^q \otimes L$, which factors as
$$
\begin{CD}
\Omega^q_{\overline{X}} \otimes L @>>>
\Omega^q_{\overline{X}}\langle D \rangle \otimes L @>>>
j_*\Omega^q_X \otimes L.
\end{CD}
$$
Since the sheaf in the middle has no cohomology in degree $p$ by
\eqref{logvan}, the composition induces $0$ on the cohomology groups
of degree $p$. We conclude that
$$
H^p(X,\Omega^q_X \otimes L) = H^0(R^p\pi_*\Omega^q_X \otimes M) = 0.
$$
Since $Y$ is affine, this proves the Theorem.
\endproof

\section{Applications.}

We will now apply Theorem~\ref{main} to derive a topological
vanishing theorem for symplectic manifolds. Assume that $X$ is a
smooth algebraic variety over $k$ equipped with a non-degenerate
closed $2$-form $\Omega \in H^0(X,\Omega^2_X)$. Assume moreover that
$X$ is equipped with a projective birational map $\pi:X \to Y$ onto
a normal algebraic variety $Y$.

\begin{lemma}\label{semismall}
The map $\pi:X \to Y$ is semismall, in other words, $\dim X \times_Y
X = \dim X$.
\end{lemma}

\proof{} For any $p \geq 0$, let $Y_p \subset Y$ be the closed
subvariety of points $y \in Y$ such that $\dim \pi^{-1}(y) \geq
p$. It suffices to prove that $\codim Y_p \geq 2p$. By \cite[Lemma
2.9]{K}, there exists an open dense subset $U \subset Y_p$ such that
the restriction $\Omega_F$ of the form $\Omega$ onto the smooth part
$F \subset \pi^{-1}(U)$ of the set-theoretic preimage $\pi^{-1}(U)
\subset X$ satisfies
$$
\Omega_F = \pi^*\Omega_U
$$
for some $2$-form $\Omega_U \in H^0(U,\Omega^2_U)$. Therefore the
rank $\rk \Omega_F$ satisfies $\rk \Omega_F \leq \dim U$. On the
other hand, since $\Omega$ is non-degenerate on $X$, we have $\rk
\Omega_F \geq \dim F - \codim F$. Together these two inequalities
give $\codim Y_p = \codim F + p \geq \dim F - \dim U + p = 2p$, as
required.
\endproof

\begin{theorem}\label{sympl}
Let $\pi:X \to Y$ be a projective birational map with smooth and
symplectic $X$. Let $y \in Y$ be a closed point, and let $E_y =
\pi^{-1}(y) \subset X$ be the set-theoretic fiber over the point
$y$. Then for odd $k$ we have $H^k(E_y,\C)=0$, while for even $k=2p$
the Hodge structure on $H^k(E_y,\C)$ is pure of weight $k$ and Hodge
type $(p,p)$.
\end{theorem}

\begin{lemma}\label{hdg}
Let $p$ be an integer, and let $V$ be an $\R$-mixed Hodge structure
with Hodge filtration $F^\hdot$ and weight filtration
$W_\idot$. Assume that $W_{2p}V = V$ and $F^pV=V$. Then $V$ is a
pure Hodge-Tate structure of weight $2p$ (in other words, every
vector $v \in V$ is of Hodge type $(p,p)$).
\end{lemma}

\proof{} Since $V = F^pV$, the same is true for all associated
graded pieces of the weight filtration on $V$. Therefore we may
assume that $V$ is pure of weight $k \leq 2p$. If $k < 2p$, we must
have $V = F^pV \cap \overline{F^pV} = 0$, which implies $V = 0$. If
$k = 2p$, the same equality gives $V = V^{p,p}$.
\endproof

\proof[Proof of Theorem~\ref{sympl}.] By Lemma~\ref{semismall},
Theorem~\ref{main} applies to $\pi:X \to Y$ and shows that
\begin{equation}\label{est}
R^p\pi_*\Omega^q_X = 0
\end{equation}
whenever $p + q > \dim X$. Since $X$ is symplectic, we have an
isomorphism $\T_X \cong \Omega^1_X$ between the tangent and the
cotangent bundle on $X$. This implies that $\Omega^q_X \cong
\Omega^{\dim X - q}_X$, and \eqref{est} also holds whenever $p > q$.

Denote by $\X$ the completion of the variety $X$ in the closed
subscheme $E_y = \pi^{-1}(y) \subset X$. Since the map $\pi:X \to Y$
is proper, by proper base change the group
$$
H^p(\X,\Omega^q_\X)
$$
for any $p$, $q$ coincides with the completion of the stalk of the
sheaf $R^p\pi_*\Omega_X^q$ at the point $y \in Y$. Therefore
$H^p(\X,\Omega^q_\X) = 0$ whenever $p > q$. The stupid filtration on
the de Rham complex $\Omega^\hdot_\X$ of the formal scheme $\X$
induces a descreasing filtration $F^\hdot$ on the de Rham cohomology
groups $H_{DR}(\X)$ which we call the weak Hodge filtration. Of
course, the associated spectral sequence does not
degenerate. Nevertheless, since $H^p(\X,\Omega^q_\X) = 0$ when $p >
q$, we have $H^k_{DR}(\X) = F^pH^k_{DR}(\X)$ whenever $k \leq 2p$.

It is well-known that the canonical restriction map
$$
H^\hdot_{DR}(\X) \to H^\hdot(E_y,\C)
$$
is an isomorphism. By definition (see \cite{D}), to obtain the Hodge
filtration on the cohomology groups $H^\hdot(E_y)$, one has to
choose a smooth simplicial resolution $E^\hdot_y$ for the variety
$E_y$ and take the usual Hodge filtration on
$H^\hdot(E^\hdot_y)$. The embedding $E_y \to \X$ gives a map
$E^\hdot_y \to \X^\hdot$, where $\X^\hdot$ is $\X$ considered as a
constant simplicial variety. The corresponding restriction map
$$
H^\hdot_{DR}(\X) \to H^\hdot_{DR}(E^\hdot_y)
$$
is also an isomorphism, and it sends the weak Hodge filtration on
the left-hand side into the usual Hodge filtration on the right-hand
side. We conclude that $H^k(E_y) = F^p(E_y)$ whenever $k \leq
2p$. It remains to recall that by definition, we have $H^k(E_y) =
W_kH^k(E_y)$, and apply Lemma~\ref{hdg}.
\endproof

To conclude the paper, we would like to note that in the particular
case when $X = T^*(G/B)$ is the Springer resolution of the nilpotent
cone $Y = \N \subset \g^*$ in the coadjoint representation $\g^*$ of
a semisimple algebraic group $G$, Theorem~\ref{sympl} has been
already proved by C. de Concini, G. Lusztig and C. Procesi in
\cite{lu}. They proceed by a direct geometric argument. As a result,
they obtain more: not only do the cohomology groups carry a Hodge
structure of Hodge-Tate type, but in fact they are spanned by
cohomology classes of algebraic cycles. This is true even for
cohomology groups with integer coefficients. Motivated by this, we
propose the following.

\begin{conj}
In the assumptions of Theorem~\ref{sympl}, the cohomology groups
$H^k(E_y,\Z)$ are trivial for odd $k$, and are spanned by cohomology
classes of algebraic cycles for even $k$.
\end{conj}

We also expect that an analogous statement holds for $l$-adic
cohomology groups, possibly even over fields of positive
characteristic.

\bigskip

\noindent
{\sc Steklov Math Institute\\
Moscow, USSR}

\bigskip

\noindent
{\em E-mail address\/}: {\tt kaledin@mccme.ru}

\end{document}